\documentstyle[12pt,amssymb]{amsart}

\newtheorem{thm}{Theorem}[section]

\newtheorem{lem}[thm]{Lemma}

\newtheorem{conj}[thm]{Conjecture}

\theoremstyle{definition}
\newtheorem{defn}{Definition}[section]

\theoremstyle{remark}
\newtheorem{rem}{Remark}[section]
\newtheorem{notation}{Notation}

\numberwithin{equation}{section}
\begin{document}

\newcommand{\thmref}[1]{Theorem~\ref{#1}}
\newcommand{\secref}[1]{\S\ref{#1}}
\newcommand{\lemref}[1]{Lemma~\ref{#1}}
\newcommand{\propref}[1]{Proposition~\ref{#1}}
\newcommand{\ban}{\!\in\!}
\newcommand{\nc}{\newcommand}
\nc{\on}{\operatorname}
\nc{\onw}{\operatornamewithlimits}
\nc{\const}{\on{const\cdot}}
\nc{\Z}{{\Bbb Z}}
\nc{\Q}{{\Bbb Q}}
\nc{\C}{{\Bbb C}}
\nc{\Pro}{{\Bbb P}}
\nc{\R}{{\Bbb R}}
\nc{\N}{{\Bbb N}}
\nc{\I}{{\cal I}}
\nc{\cL}{{\cal L}}
\nc{\J}{{\cal J}}
\nc{\lpg}{\widehat{LG}}
\nc{\A}{\operatorname{A\cdot\!}}
\nc{\al}{\alpha}
\nc{\De}{\Delta}
\nc{\eps}{\varepsilon}
\nc{\io}{\iota}
\nc{\La}{\Lambda}
\nc{\la}{\lambda}
\nc{\th}{\theta}
\nc{\khi}{\chi}
\nc{\om}{\omega}
\nc{\Om}{\Omega}
\nc{\isom}{\cong}
\nc{\bs}{\bar{s}}
\nc{\lar}{\longrightarrow}
\nc{\arr}{\rightarrow}
\nc{\SU}{\operatorname{SU}}
\nc{\F}{{\bf F}}
\nc{\GL}{\on{GL}}
\nc{\SL}{\on{SL}}
\nc{\U}{\on{U}}
\nc{\SO}{\on{SO}}
\nc{\Spin}{\on{Spin}}
\nc{\M}{{\goth M}}
\nc{\Sym}{\on{Sym}}
\nc{\ahat}{\hat{A}}
\nc{\spec}{\on{pecm}}
\nc{\End}{\on{End}}
\nc{\fel}{\frac{1}{2}}
\nc{\maps}{\longmapsto}
\nc{\zi}{\frac{1}{z}}
\nc{\Res}{\onw{Res}}
\nc{\loc}{\on{loc}}
\nc{\Ext}{\on{Ext}}
\nc{\Hom}{\on{Hom}}
\nc{\Pic}{\on{Pic}}
\nc{\cH}{\on{H}}
\nc{\Ch}{\on{Ch}}
\nc{\Todd}{\on{Todd}}
\nc{\res}{\on{Res}}
\nc{\Spec}{\on{Spec}}
\nc{\Wr}{W_r}
\nc{\perm}{\on{perm}}
\nc{\signs}{\on{signs}}
\nc{\Wp}{W_r^{\perm}}
\nc{\Ws}{W_r^{\signs}}
\nc{\Gr}{\on{Gr}}
\nc{\hook}{\hookrightarrow}
\nc{\frwd}{\fracwithdelims()}
\nc{\be}{\begin{equation}}
\nc{\Vm}{\on{Vm}}
\nc{\Tr}{\on{Tr}}
\nc{\PO}{\on{PO}}
\nc{\tPO}{\tilde{\PO}}
\nc{\tPG}{\tilde{\PG}}
\nc{\PG}{\on{PG}}
\nc{\tensor}{\otimes}
\nc{\rank}{\on{rank}}
\nc{\ds}{\displaystyle}
\nc{\sL}{{\cal L}}
\nc{\sE}{{\cal E}}
\nc{\sS}{{\cal S}}
\nc{\sT}{{\cal T}}
\nc{\cO}{{\cal O}}
\nc{\tT}{\tilde T}
\nc{\zin}{z^{-1}}
\nc{\pone}{\Pro^1}
\nc{\slope}{\on{sl}}
\nc{\gb}{\bar{g}}
\nc{\g}{{\goth g}}
\nc{\del}{\delta}
\nc{\gs}{{\goth g}^*}
\nc{\gt}{{\goth t}}
\nc{\ts}{\gt^*}
\nc{\W}{\Lambda^*}
\nc{\mw}{\mu_W}
\nc{\cl}{\chi_\lambda}
\nc{\ar}{\longrightarrow}
\nc{\eSp}{\operatorname{Sp}}
\nc{\Gtwo}{\operatorname{G}_2}
\nc{\Ffour}{\operatorname{F}_4}
\nc{\Wt}{\cal W}
\nc{\goC}{{\goth C}}
\nc{\<}{\langle}
\nc{\?}{\rangle}
\nc{\alc}{{\goth a}}
\nc{\K}{\on{K}}
\nc{\indu}{\on{ind}}

\title[Verlinde formulas]{The combinatorics of the Verlinde formulas}

\author{Andr\'as Szenes}
\address{Massachusetts Institute of Technology, Department of Mathematics,
Cambridge, MA 02139}
\thanks{Research was partially supported by an NSF grant.}

\email{szenes@@math.mit.edu}

\maketitle

\section{Introduction}

In this short note we discuss the origin and  properties of the Verlinde
formulas and their connection with the intersection numbers of moduli
spaces. Given a simple, simply connected Lie group $G$, the Verlinde
formula is an expression $V^G_k(g)$ associated to this group depending
on two integers $k$ and $g$.  For $G=\SL_2$ the formula is  \be
\label{evenintro} V_k^{\SL_2}(g) = \sum_{j=1}^{k-1}
\frwd{k}{2\sin^2\frac{j \pi}{k}}^{g-1}. \end{equation}

We describe $V^G_k$ for general groups in \secref{fusion}. These formulas
were first written down by E. Verlinde \cite{fusion} in the context of
Conformal Field theory. The interest towards them in Algebraic Geometry
stems form the fact that they give the Hilbert function of moduli spaces of
principle bundles over projective curves.  More precisely, let $C$ be a
smooth projective curve of genus $g$, and let $\M_C^G$ be the moduli space
of principal $G$-bundles over $C$ (cf. e.g. \cite{kunara} and references
therein). Then there is an ample line bundle $\cL$ over $\M_C^G$ such that
\be \label{main}
\dim \cH^0(\M_C^G,\cL^k) = V_{k+h}^G(g),
\end{equation}
where $h$ is the dual Coxeter number of $G$. This statement requires
some modifications for a general   simple $G$, but it holds for $\SL_n$
(\cite{BSz,faltings,bela,kunara}).

Proving \eqref{main} is important, but in this paper we will
address a different question: what can be said about the moduli spaces
knowing \eqref{main}?  Accordingly, first we concentrate on
understanding the formula.

Two rather trivial aspects of \eqref{main} are that
\begin{itemize}
\item  $V^G_k(g)$ is integer valued,
\item $V^G_k(g)$ is a polynomial in $k$.
\end{itemize}
Note that looking at the formula itself, none of this is
obvious. Our goal  is to explain these properties and
connect them to the intersection theory of $\M^G_C$.

The paper is structured as follows: in \secref{fusion} we discuss some of
the ideas of Topological Field Theory, which explain the structure of the
formula for general $G$ and show its integrality (cf.
\cite{kac,phd,gepner,arnb}). In \secref{resi} we give our main result, a
residue formula for $V^G_k$ for $G=\SL_n$. Such a formula gives an explicit
way of calculating the coefficients of $V^G_k$ as a polynomial in $k$.
Finally, in
\secref{application} we give an application of our formulas: a ``one-line
proof of'' \eqref{main}.

This paper is intended as an announcement and overview. As a result, few
proofs will be given, and even most of those will be sketchy. A more
complete treatment will appear separately.

{\bf Acknowledgements.} I am grateful to Raoul Bott, my thesis advisor,
for suggesting to me this circle of problems and helping me with advice
and ideas along the way. I would like to thank Noam Elkies for useful
discussions.

I am thankful to the organizers of the Durham Symposium on Vector bundles,
Peter Newstead and Bill Oxbury, for their help and for the opportunity to
present my work.

\section{Topological Field Theory and Fusion Algebras} \label{fusion}

This section is independent from the rest of the paper. It contains a
quick and rather formal overview of the structure of Topological Field
Theories \cite{segal,atiyah} and Verlinde's calculus \cite{fusion}.

Consider a finite dimensional vector space $F$ (the space of fields)
with a marked element $1\in F$ (the vacuum). Assume that a  number
$F(g)_{v_1,v_2,\dots v_n}$ (correlation functions) is associated to every
topological Riemann surface  of genus $g$, with elements of the algebra
$v_1,v_2,\dots v_n\in F$ inserted at $n$ punctures, which satisfies the
following axioms:  \begin{description} \item[Normalization]
$F(0)_{1,1,1} = 1$, \item[Invariance] $F(g)_{v_1,\dots} =
F(g)_{1,v_1,\dots}$, \item[Linearity] $F(g)_{v_1,\dots}$ is linear in
$v_i$. \end{description} Introduce the symmetric linear 3-form
$\om:F\tensor F\tensor F\arr \C$ by $\om(u,v,w) = F(0)_{u,v,w}$, the
bilinear form $(u,v)=F(0)_{u,v}$ and the trace $\int u = F(0)_{u}$.
Assume that $(,)$ is {\bf non-degenerate}, and fix  a pair of  bases
$\{u_i,u^i\}$ of $F$, dual with respect to this form, that is
$(u_i,u^j)=\delta_{ij}$.

\begin{description}
\item[Verlinde's fusion rule]  $F(g)_{v_1,\dots} = \sum_{i}
F(g-1)_{u_i,u^i,v_1\dots}. $
\end{description}

One can extend $F$ to disconnected surfaces by  the  axiom:
\begin{description}
\item[Multiplicativity] $F$ is
multiplicative under disjoint union.
\end{description}

\begin{rem} These axioms serve as an algebraic model of certain
relations among the Hilbert functions of various moduli spaces. The number
$F(g)_{v_1,v_2,\dots v_n}$ represents the dimension of the space of
sections of a certain line bundle over a moduli space of parabolic bundles
with weights depending on the insertions $v_1,v_2,\dots v_n$. The fusion
axiom describes how the space of sections of a line bundle decomposes over
a family of curves degenerating to a nodal curve (cf.
\cite{faltings,tsueya}).
\end{rem}

\begin{lem} \label{twisted:moments}  The axioms above define the
structure of an associative and commutative algebra on $F$, by the
formula $vw=\sum_i\om(v,w,u^i)u_i$, compatible with $(,)$ and $\int$.
Then if we denote the invariantly defined element $\sum_i u_i u^i\in F$
by $\al$, we have
\be \label{numberofCB} F(g)_{v_1,v_2,\dots,v_n} = \int
\al^g v_1v_2\dots v_n. \end{equation} \end{lem}

Now assume in addition that the algebra $F$ is {\bf semisimple}. Then it
has the form $ F\cong L^2(S,\mu),$ where $S=\Spec F$ is a finite set
and the complex measure $\mu$ can be given via a function $\mu:
S\arr\C$.

The elements of $F$ become functions on $S$ and the trace
$\int$ turns out to be the actual integral with respect to
$\mu$. Now take the following pair of dual bases:
$\{\del_s,\del_s/\mu(s)|\; s\in S\}$, where $\del_s(x)=\del_{sx}$. We
call this the spectral basis. Using this basis and \eqref{numberofCB},
we obtain the following formula:
\be \label{diagonal:formula}
 F(g) = \sum_{s\in S} \mu(s)^{1-g}.
\end{equation}

This formula resembles \eqref{evenintro}, but what is the appropriate
algebra?

\subsection{Fusion algebras}

Here we construct the fusion algebras for arbitrary simple, simply
connected Lie groups. First we need to introduce some standard
notation.

\begin{notation} In this paragraph we will use the compact form of
simple Lie groups, still denoting them by the same letter.  Thus let $G$ be
a compact, simply connected, simple Lie group, $\g$ its complexified Lie
algebra, $T$ a fixed maximal torus, and $\gt$ the complexified Lie algebra
of $T$.  Denote by $\La$ the unit lattice in $\gt$ and by $\Wt\in\ts$ its
dual over $\Z$, the weight lattice.  Let $\De\in\Wt$ be the set of roots
and $W$ the Weyl group of $G$. A fundamental domain for the natural action
of the Weyl group on $T$ is called an alcove; a fundamental domain for the
associated action on $\ts$ is called a chamber. We will use the
multiplicative notation for weights and roots, and think of them as
characters of $T$.  The element of $\ts$ corresponding to a weight $\la$
under the exponential map will be denoted by $L_\la$.

Fix a dominant chamber $\goC$ in $\ts$ or a corresponding alcove $\alc$
in $T$. This choice induces a splitting of the roots into positive
($\De^+$) and negative ($\De^-$) ones. For a weight $\la$, denote its
Weyl antisymmetrization by $\A\la=\sum_{w\ban W} \sigma(w)w\cdot\la$,
where $\sigma:W\arr\pm1$ is the standard character of $W$. According to
the Weyl character formula, for a dominant weight $\la$, the character
of the corresponding irreducible highest weight representation is
$\chi_\la=\A\la\rho/\A\rho$, where $\rho$ is the square root of the
product of the positive roots.

The ring $R(G)$, the representation ring of $G$, can be identified with
$R(T)^W$, the ring of Weyl invariant linear combinations of the weights.
Denote by $d\mu_T$ be the normalized Haar measure on $T$.  If we endow
$T/W$ with the Weyl measure $$d\mw = \A\rho \A\bar\rho \;d\mu_T,$$
then $R(G)$ becomes a pre-Hilbert space with orthonormal basis
$\{\cl\}$, i.e. one has \newline
$\int_{T/W}\cl\chi_\mu\,d\mw = \delta_{\la,\bar{\mu}}.$
 \qed \end{notation}

We need to introduce an integer parameter denoted by $k$ called
the {\em level}, which can be thought of as an element of
$\cH^3(G,\Z)\cong\cH^4(BG,\Z)\cong\Z$, and in turn can be identified
with a Weyl-invariant integral inner product on $\gt$.

The {\em basic} invariant inner product on $\gt$ corresponding to
$k=1$ is specified by the condition $(H_\th,H_\th)=2$, where $H_\th\in\gt$
is the coroot of the highest root $L_\th$. It has the following
properties (see \cite[\S 6]{kac},\cite[Ch.4]{loopgr}):
\begin{itemize}
\item For the induced inner product on $\ts$, we have $(L_\th,L_\th)=2$.
\item For $\la\ban\Wt$,
the inner product $(L_\th,L_\la)$ is an integer,  and (,) is the smallest
inner product with this property.
\item The Killing form is equal to $-2h (,)$, where $h = (L_\th,L_\rho)+1$ is
the dual Coxeter number of $G$.
\end{itemize}
The basic inner product also gives an identification $\nu:\ts\arr\gt$
between $\ts$ and $\gt$, by the formula $\beta(x)=(\nu(\beta),x)$.

\subsection{The  simply-laced subgroup}

Let $\De_l\in\De$ be the set of long roots of $G$.  Denote by $\Wt_r$ the
lattice in $\ts$ generated by $\De$, and by $\Wt_l$ the lattice generated
by $\De_l$. By definition $\La$ is the dual of $\Wt$ over $\Z$ with respect
to the canonical pairing $\<\, ,\?$ between $\ts$ and $\gt$. The dual
of $\Wt_r$ is the  center lattice in $\gt$. Denote
the dual of $\Wt_l$ by $\La_l$.

The root system $\De_l$ corresponds to a subgroup $G_l$ of $G$ with
maximal torus $T$ and Weyl group $W_l\subset W$, which is generated by
reflections corresponding to the elements of $\De_l$.  Denote the center
of $G_l$ by $Z_l$. Then $Z_l$ can be described as the set of elements of $T$
invariant under $W_l$, and we have $\exp^{-1}\La_l=Z_l$. It is important
to note that in view of the second property of $(,)$ above, $\nu:\Wt_l \maps
\La$ is an isomorphism.
Since $\Wt$ is paired to $\La$ and $\Wt_l$ is paired to
$\La_l$ over $\Z$, it follows that $\nu:\Wt\maps
\La_l$ is also an isomorphism. Then the map
$\exp\cdot\nu:\{\al\ban\goC|\,(L_\th,\al)\leq 1\}\maps\alc$ is a
bijection.

Naturally, if $G$ is simply laced, then $G_l=G$. For the non-simply laced
groups one has the following subgroups:
\begin{itemize}
\item $\Spin_{2n}\subset\Spin_{2n+1} $
\item $\SU_2^n\subset\eSp_n $
\item $\SU_3\subset\Gtwo $
\item $\Spin_8\subset\Ffour $
\end{itemize}

\subsection{The definition of the fusion algebra}

We give a different definition from the standard one via co-invariants
of infinite dimensional Lie algebras \cite{tsueya}, but one which is
very natural from the point of view of representation theory.

To motivate the construction,  recall the procedure of holomorphic
induction \cite{bott}: the flag variety $\F=G/T$ has a complex structure
and every character $\la$ of $T$ induces a holomorphic equivariant line
bundle $\cL_\la=G\times_\la \C$ over $\F$. Then one can define the
induction map $\I:R(T)\arr R(G)$ as a homomorphism of additive groups by
the formula $\la\mapsto \sum (-1)^{i}\cH^i(F,\cL_\la)$, where the
cohomology groups in the latter expression are thought of as
$G$-modules. The Borel-Weil-Bott theorem then says that
\begin{align}
&\text{for }\la\text{ dominant }\I(\la)=\chi_\la\\
\intertext{and}
&\I(\la') = \sigma(w)\I(\la), \text{ whenever } \la'\rho =
w(\la\rho)\text{ for some } w\in W. \label{indrel}
\end{align}
Note that $\I$ is not expected to be a ring homomorphism.

This procedure applies to the loop group $\lpg$ as well
(\cite{loopgr,kumar,mathi}). Once the action of the central elements is
fixed as $c\mapsto c^k$, where $k\in \N$ is the level, again, we have a map
$\bar{\I}:R(T)\mapsto R_k(\lpg)$. This last object $R_k(\lpg)$ has only
additive structure, since the tensor product of two level $k$
representations has level $2k$. The role of the Weyl group is played by the
affine Weyl group $W_k$ obtained by adjoining to $W$ the translation by
$(k+h) L_\th$. Again the Borel-Weil-Bott theorem applies, and
\eqref{indrel}, with $W$ replaced with $W_k$, gives a description of the
kernel of $\bar{\I}$.

Since $W\subset W_k$, the map $\bar{\I}$ factors through $\I$, and as a
result we have a map $\J:R(G) \mapsto R_k(\lpg)$. It is easy to see from
the \eqref{indrel} that the set of characters: $\Xi _k = \{\cl |\,
(L_\th,L_\la)\leq k\}$ forms a basis of $R_k(\lpg)$ if we identify $\J(\cl)$
with $\cl$.

\begin{lem} \label{quot}
The additive group $R_k(\lpg)$ can be endowed with a ring
structure $F_k^G$ so that the map $\J$ becomes a homomorphism of rings.
\end{lem}

The algebra $F_k^G$ is called the {\em fusion algebra} of $G$ of level $k$.
As noted above, we can consider $\Xi_k$ to be a basis of $F^G_k$. Endow
$F_k^G$ with the trace function $\int$ by the formula $\int\cl=0$ except
for the trivial character $\chi_1$, which has trace equal to 1. Also
note that since $\Spec(R(G))=T/W$ and $F_k^G$ is a quotient of $R(G)$, we
expect $\Spec(F_k^G)\subset T/W$.

\begin{lem} \label{measure}
$F_k^G$ can be identified with ``$L^2$'' of the finite normalized measure
space $Z_k=\{t\ban T |\, t^{k+h}\ban Z_l, \, t\;\text{\rm is
regular}\}/W$, with measure $d\mu_k$ given by the function
\be \label{lu}
\dfrac{\A\rho(t)\A\bar\rho(t)}{|Z_l|(k+h)^r}. \end{equation}
\end{lem}

Note the surprising fact that the discrete measure remains unchanged up
to a normalization factor as $k$ varies.

Now we can define the  quantity $V^G_k(g)$ which appeared in
\eqref{main} as the number associated to a Riemann surface of
genus $g$ and the fusion algebra $F^G_{k-h}$. Combining
\eqref{diagonal:formula} and \eqref{lu} we obtain
$$ V^G_k(g) = \sum_{t\in T_r/W, t^k\in Z_l}
\frwd{|Z_l|k^r}{\A\rho(t)\A\bar\rho(t)}^{g-1},$$
where $T_r$ is the set of regular elements of $T$. This can be easily
seen to give \eqref{evenintro} for the case of $G=\SU_2$. Indeed,
embedding the maximal torus of $\SU_2$ into $\C$ as the unit circle,
we have: $\rho(z)=z$, $Z_l=\pm 1$, $h=2$ and $A(z) = 1/z$.

Finally, note that since the relations in the fusion algebras given in
\eqref{indrel} have integer coefficients, and the $\{\khi_\la\}$ form an
orthonormal basis of $F$, we see that $\int\al^g$ from \eqref{numberofCB}
has to be an integer. This proves that $V^G_k(g)$ is an integer for all
groups and values of $k$ and $g$.

\begin{rem} That this  definition of the fusion algebras is equivalent
to the standard one via coinvariants of current algebras
\cite{tsueya} can be shown to be
equivalent to Verlinde's conjecture on the diagonalization of the fusion
rules.  which gives a formula for the product in $F^G_k$ using the
$S$-matrix. The definition given above is simpler to use for calculations
and it gives the correct prescription for non simply connected groups (see
also \cite{gepner}).
\end{rem}

\section{Residue formulas} \label{resi}

 In this section we study $V_k^G(g)$ as a function of $k$. We show that
$V_k^G(g)$ is a polynomial in $k$ for $G=\SL_3$, and  give a simple
formula for the coefficients of this polynomial. The generalization of
these results to $\SL_n$ is  straightforward.

Consider the case $G=\SL_2$ first. Again, as at the end of the previous
section embed the maximal torus of $\SL_2$ into $\C\subset\pone$.

Consider the differential form $$\mu=\frac{dz}{z}\frac{z+\zin}{z-\zin}
\text{ on }\Pro^1.$$ This form has simple poles: at $z=\pm 1$ with
residue 1, and at $z=0,\infty$ with residue $-1$. Thus if we pull back
$\mu$ by the $k$-th power map we obtain a differential form $\mu_k$
with poles at the $2k$-th roots of unity and residues $+1$, and simple
poles at $z=0,\infty$ with residue $-k$ . It is given by the following
formula: $$ \mu_k = k
\frac{dz}{z}\frac{ z^{k}+z^{-k}}{z^{k}-z^{-k}}.$$
Note that $\mu_k$ is invariant under multiplication by a $2k$th root
of unity and under the Weyl reflection  $z\arr 1/z$.

Now suppose we have a function $f(z)$, with poles only at $z=\pm 1$,
vanishing at 0 and $\infty$, and invariant under the substitution
$z\arr z^{-1}$. Then by applying the Residue Theorem to
the differential form $f\mu_k $ and using the Weyl symmetry at hand,
we have
$$\sum_{j=1}^{k-1}f(\exp(\pi Ij/k)) =
-\Res_{z=1} \mu_k f(z),$$
where $I^2=-1$.
Applying this argument to the function
$$f(z)=\frwd{2k}{-(z-z^{-1})^2}^{g-1}$$ we obtain the formula
$$V^{\SL_2}_k(g) =(-1)^{g} (2k)^{g-1}
\Res_{z=1} \frac{k\, dz}{z} \frac{z^{k}+z^{-k}}{z^{k}-z^{-k}}
\frwd{1}{(z-z^{-1})^2}^{g-1}.$$

Now using the invariance of the residue under substitutions we can
obtain different formulas for $V^{\SL_2}_k(g)$. For example, the
polynomial nature of $V^{\SL_2}_k(g)$ becomes transparent if we
perform the substitution $z\arr \exp(Ix)$: $$V^{\SL_2}_k(g) =
-(2k)^{g-1}\Res_{x=0} \frac{k\cot(kx)\,dx}{(2\sin x)^{2(g-1)}}.$$   It
is easy to check using this formula that the degree of $V^{\SL_2}_k(g)$
as a polynomial in $k$ is $3(g-1)$, which, as expected, coincides with
the dimension of $\M^{\SL_2}_C$.

Before we proceed, we need an understanding of higher dimensional
residues. The notion that a top dimensional differential form has an
invariantly defined number assigned to it, does not carry over to the
higher dimensions. The correct object in $\C^n$ is
$\Res:\cH^n_{\loc}(\Om^n,\C^n) \mapsto \C$ mapping from the $n$th local
\v Cech cohomology group in a neighborhood of 0 with values in  holomorphic
$n$-forms  to
complex numbers. To define this map let $\om$ be a meromorphic $n$-form
defined in a neighborhood of 0 in $\C^n$. Then $\om$ can be represented in
the form $dz_1\,dz_2\,\dots dz_n h(z)/f(z)$ where $f$ and $h$ are
holomorphic functions. The additional data necessary to represent an
element of $\cH^n_{\loc}(\Om^n,\C^n)$ is a splitting of $f$ into the
product of $n$ functions $f=a_1a_2\dots a_n$. Such a splitting defines $n$
open sets in the a neighborhood of 0: $A_i = \{a_i\neq 0\}$.  These define
a local \v Cech cocycle. A detailed explanation of this and an algorithm to
calculate the residue can be found in \cite{grha,harts}.

We will call a differential $n$-form with such a splitting a {\em
residue form}.

\begin{defn} A non-trivial residue form $\om$ is called flaglike if $a_i$ only
depends on $z_1,\dots z_i$. This notion depends on choice and the order of
the coordinates $z_1,\dots, z_n$.  \end{defn}

\begin{lem} \label{flag}
 Let $\om$ be a flaglike residue form. Then
$$ \Res(\om) = \Res_{z_n}\dots\Res_{z_1} \om.$$
Here $\Res_{z_i}$ is the ordinary 1-dimensional residue, taken assuming
all the other variables to be constants.
\end{lem}
The proof is  straightforward. Note that the order of the variables is
important, while there is some freedom in the way the denominator is split
up.

For simplicity we restrict ourselves to the case of $\SL_3$. According
to \lemref{measure} and \eqref{diagonal:formula} the Verlinde formula
can be written as
\be \label{eseltri}
 V_k^{\SL_3}(g) = (3k^2)^{g-1}\sum_{i,j,k-i-j > 0}
(8\sin(i\pi/k)\sin(j\pi/k)\sin((i+j)\pi/k))^{-2(g-1)}
\end{equation}

Now we can write down the main result of the paper:
\begin{thm} \label{result}
\begin{multline} \label{fores}
 V_k^{\SL_3}(g) = (3k^2)^{g-1} \Res_{Y=1} \Res_{X=1}
\frac{X^k+X^{-k}}{X^k-X^{-k}}\;\frac{Y^k+Y^{-k}}{Y^k-Y^{-k}}\times\\
\times\frac{(-1)^{g-1}}{((X-X^{-1})(Y-Y^{-1})(XY-(XY)^{-1}))^{2(g-1)}}
\frac{k^2\,dX\,dY}{XY}.
\end{multline}
\end{thm}

The proof is analogous to the case of $\SL_2$. Denote the residue form in
\eqref{fores} by $\om_k(g)$.
Again at the points
$p_{ij} = (e^{iI\pi/k},e^{jI\pi/k})$, with $i,j,k-i-j\ge 1$ the
residues of $\om_k(g)$ reproduce the sum \eqref{eseltri}. However,
 now it is not immediately obvious that the residue theorem can
localize this sum at the point $(1,1)$, since the residue form in
\eqref{fores} has non-trivial residues at other points as well. To
illustrate the situation consider the matrix $M_k$ whose $(i,j)$th entry
is the residue of  $\om_k(g)$ taken at the point $p_{ij}$
instead of $(1,1)$, where $i,j=0,\dots,k-1$.

Example for $g=2$: $$ M_6 =
\left [\begin {array}{cccccc} 166&-45&-29&-18&-29&-45
\\\noalign{\medskip}-45&36&9&9&36&-45\\\noalign{\medskip}-29&9&4&9&-29
&36\\\noalign{\medskip}-18&9&9&-18&9&9\\\noalign{\medskip}-29&36&-29&9
&4&9\\\noalign{\medskip}-45&-45&36&9&9&36\end {array}\right ]
. $$

 We can apply the Residue Theorem to ``each column'' by fixing a value
of $X$.  By degree count, one can see that $\om(g)$ has trivial residues at
$Y=0,\infty$ and this implies
that the sum of the entries in each column of $M_k(g)$ is 0. Next,  note that
$M_k(0,i)=M_k(i,0)$, since these residues are {\em split}, i.e. they have
the form $dX\,dY\, X^{-m}Y^{-n} f(X,Y)$, where $f$ is holomorphic at the
point where the residue is taken.

Now to prove the Theorem it is sufficient to show that
$M_k(j,0)=M_k(j,k-j)$ for every $j>0$.  To see this, note that both residues
are simple (first order) in $X$ at $\al=\exp(j\pi/k)$. This means that
after taking the $X$-residue, we are left with the form $$ \om_\al = \const
\frac{dY}{Y}
\frac{Y^k+Y^{-k}}{Y^k-Y^{-k}} \frac{1}{((Y-Y^{-1})(\al
Y-\al^{-1}Y^{-1}))^{2(g-1)}}.$$ The two numbers we need to compare are the
residues of this form at $0$ and $\al$ respectively. But these two
residues clearly coincide since $\om_\al$ is invariant under the
substitution $Y\arr \al^{-1}Y^{-1}$. \qed

The formula for $G=\SL_n$ reads as follows:
\be \label{slnverl}
V^G_k(g) = (-1)^{n-1+(g-1)|\Delta^+|} (nk^{n-1})^{g-1} \Res_{X_{n-1}=1}\dots
\Res_{X_1=1}
W^{-2(g-1)}
\prod_{i=1}^{n-1}  \frac{X_i^k+1}{X_i^k-1} \cdot \frac{k\,dX_i}{2X_i}
\end{equation}
where $X_i,\; i=1,\dots,n-1$ are the simple (multiplicative) roots and $W =
\prod_{\al\in\Delta^+} (\al^{\fel}-\al^{-\fel})$ is the Weyl
measure.

As we pointed out after \lemref{flag}, the ordering of the variables
matters when taking the subsequent residues. In the special case of
$G=\SL_3$ this ordering does not matter (i.e. $M_k(i,j)=M_k(j,i)$), but for
higher rank groups a finer argument is necessary.

Finally, note that similarly to the case of $\SL_2$, \eqref{slnverl} gives
a simple prescription for calculating the coefficients of $V^G_k(g)$ as a
polynomial in $k$, via the exponential substitution. For example, for
$\SL_3$ we obtain
\be \label{verlslt} V^{\SL_3}_k(g)= (3k^2)^{g-1} \Res_{y=0}\Res_{x=0}
\frac{k^2cot(kx)\cot(ky)\, dx\,dy}{(8\sin(x)\sin(y)\sin(x+y))^{2(g-1)}}.
\end{equation}

A different generating function was obtained for the case of $G=\SL_3$ by
Zagier \cite{zagverl}.

\section{Multiple  $\zeta$-values and intersection numbers of moduli spaces}
\label{application}

In this final section we show how \eqref{main} and \eqref{fores} can be
related via the Riemann-Roch formula to Witten's conjectures on the
intersection numbers of moduli spaces. Our argument below gives a quick
proof of \eqref{main} for $\SL_n$ assuming Witten's formulas. This is a
generalization of the work of Thaddeus who considered the case of $\SL_2$
\cite{thad}.

\subsection{Multiple $\zeta$-values and intersection numbers of moduli
spaces}

Consider the case of $\SL_2$ first. If we want to find the asymptotic
behavior of $V^{\SL_2}_k$ for large $k$, the best way to think about the
formula is that it is a discrete approximation to the (divergent)
integral $\int_0^1 \sin(\pi x)^{-2(g-1)}\, dx$. To find the leading
asymptotics, we can replace $\sin(x)$ by $x$, and taking the large $k$
limit we obtain: $ (k/2)^{g-1}\sum_{j=1}^\infty (k/(j\pi))^{2(g-1)} =
k^{3(g-1)}\zeta(2(g-1))/(2^{g-1}\pi^{2(g-1)})$. This can be easily
proven, and in fact, a generalization of  this formula  for arbitrary
groups appeared in Witten's work \cite{wittzeta}.

Below we will concentrate on the case of $\SL_3$, however the formulas can
be extended to $\SL_n$ as well.

If we perform the  trick above for $\SL_3$, up to a constant, the leading
behavior of the Verlinde formula appears to be  $$ V^{\SL_3}_g(k) \sim
\const k^{8(g-1)}/\pi^{6(g-1)} \sum_{i,j=1}^\infty  (ij(i+j))^{-2(g-1)}.
$$

One can write down more general sums, e.g.:
$$S(a,b,c) =  \sum_{i,j=1}^\infty i^{-a} j^{-b} (i+j)^{-c},$$
closely related to the so-called {\em multiple zeta values} \cite{zagzet}.

It was discovered by Witten that all intersection numbers of moduli
spaces are given by combinations of multiple $\zeta$-values
\cite{wittvol,wittzeta}. Below we give a couple of useful formulas for them.
We restrict ourselves to the case $S(2g,2g,2g)$ for simplicity. Similar
formulas exist in greater generality.

\begin{lem} \mbox{}
\be \label{bern}
S(2g,2g,2g) = \fel\int_0^1 \bar{B}_{2g}(x)^3 \, dx,
\end{equation}
 where $\bar{B}_n(x)$ is a modified {\rm n}th Bernoulli polynomial,
$\bar{B}_n(x)= -(2\pi I)^n B_n(x)/n!$.
\be \label{zetares}
S(2g,2g,2g) =
\frac{1}{3}\Res_{(0,0)} \cot(x)\cot(y)  (xy(x+y))^{-2g}.
\end{equation}
\end{lem}

Sketch of Proof: The first formula follows from the definition of the
Bernoulli polynomials:
$$ \bar{B}_n(x) = \sum_{j\neq 0} e^{2Ij\pi x}/j^n. $$
Indeed, substituting this into $\fel\int_0^1 B_{2g}(x)^3 \, dx,$ one obtains
$S(2g,2g,2g)$ on the nose; the coefficient $1/2$ is a combinatorial
factor.

The proof of the second formula is similar to the proof of
\thmref{result}. On has to apply the Residue Theorem in two steps. That
the residue at infinity vanishes follows from the  expansion $\cot(x)$:
$$\pi\cot(\pi x) = \sum_{n\in\Z} (x-n)^{-1}.\qed $$

\subsection{Intersection numbers of the moduli spaces}  First we  recall
some facts about the cohomology of the moduli spaces. We will ignore that
the moduli spaces are not smooth in general, and accordingly, we will
assume the existence of a universal bundle, Riemann-Roch formula, etc.
However, formulas analogous to \eqref{result} exist for the smooth moduli
spaces as well (e.g. when the degree and rank are coprime for $\SL_n$), and
all of our statements are rigorous for these cases. Some of the singular
moduli spaces (e.g. vector bundles) can be handled using the methods of
\cite{BSz}. We will also ignore certain difficulties which arise for
$\Spin_{n}$, $n>6$, and the exceptional groups, where the ample line
bundle exists only for $k=0\mod l$, for some $l$, depending on the type of
the group. In these cases the Verlinde formula is a polynomial only when
restricted to these values. Thus what follows should perceived as a
scheme of a proof, which works as it is in some cases, but requires
modification and more work in greater generality.

There is a universal principal $G$-bundle $U$ over the space
$C\times\M^G_C$, which induces a map $\M^G_C\arr BG$, and consequently
a map $s:\cH^*(BG)\arr \cH^*(\M^G_C)\tensor\cH^*(C)$.

Recall that $\cH^*(BG) \cong\Sym(\g^*)^G$, the space of $G$-invariant
polynomial functions on $\g$. This is a polynomial ring itself in
$\rank(G)$ generators and it is isomorphic by restriction to
$S^G=\Sym(\gt^*)^W$. For every $\al\in\cH_i(C)$ and $P\in S^G$ we obtain a
cohomology class of $\al\cap s(P)\in\cH^{2i-j}(\M^G_C)$, the
$\al$-component of $s(P)$.  In fact, $s$ induces a map $\bs:\cH_*(C)\cap
S^G\arr\cH^*(\M^G_C)$, where $\cH^*(C)\cap S^G$ is the free commutative
differential algebra generated by the ring $S^G$ and the differentials
of negative degree modeled on $\cH_*(C)$. For the case of $\SL_n$ and
coprime degree and rank it is known that $\bs$ is surjective \cite{AB,kirwan}.
In
particular, denoting the fundamental class of $C$ by $\eta_C$, and the
basic invariant scalar product from \secref{fusion} by $P_2$ we obtain a
class $\om = \eta_C\cap s(P_2)\in \cH^2(\M^G_C),$ which turns out to be the
first Chern class of the line bundle from \eqref{main}. To simplify the
notation, below we omit the map $s$ and also $\al$ if $\al=1$, when writing
down the classes $\cH^*(\M^G_C)$. Thus $1\cap s(P)$ will be denoted simply
by $P$.

Any power series in the variables $\al\cap P$ can be integrated over
$\M^G_C$ and these numbers are called the intersection numbers of the
moduli space. Naturally, only the terms of degree $\dim \M^G_C=
\dim(G)(\on{genus}(C)-1)$ will contribute.

Witten, using non-rigorous methods, gave a complete description of these
intersection numbers in the most general case \cite{wittzeta}.  His
formulas are  combinations of multiple $\zeta$-values, and are rather difficult
to calculate. In this paper, we will focus only on a subset of these
intersection numbers, which are of the form $\int_{\M} \om^l P$, where
$l\in\N$ and $P$ is a not necessarily homogeneous Weyl-symmetric function
on $\gt$.

\begin{conj} \label{conjecture}
For every group $G$, there exists a residue form $\Omega^G$ depending on
$g$, defined in a neighborhood of $0\in\gt$, the Cartan subalgebra of $G$,
such that
\be \label{conj}
\int_{\M} e^{\om} P = \Res_{\on{at}\; 0\in\gt} \Omega^G P,
\end{equation}
\end{conj}

For $G=\SL_n$,
\be
\label{intersect}
\Omega = n^{g-1} \Res_{x_{n-1}=0}\dots\Res_{x_1=0} \prod_{\al\in
\Delta^+} L_\al^{2(g-1)} \prod_{i=1}^{n-1} \cot(x_i) \, dx_i,
\end{equation}
where the $x_i$-s are halves of the simple (additive) roots of $\SL_n$,
ordered according to the Dynkin diagram.

Let us write down the formula for $\SL_3$ more explicitly and inserting the
``grading'':
\be \label{special}
\int_{\M} e^{k\om} P = (3k^2)^{g-1}
\Res_{y=0} \Res_{x=0}
\frac{k^2\cot(kx)\cot(ky)\,dx\,dy}{(8xy(x+y))^{2(g-1)}}P
\end{equation}

\begin{rem}
It can be shown that \eqref{intersect} is consistent with  Witten's
formulas. We will not give the proof here, but note that the link
between the two types of formulas is given by equalities like
\eqref{zetares}.

At the moment we do
not know $\Omega^G$ for general $G$.

Our formulas seem to  be related to those given in the
works of Jeffrey and Kirwan \cite{jeffkir,jeffkirnew}.  \qed
\end{rem}

Finally, we present another evidence for \eqref{intersect}:
the consistency with the Verlinde formula. First we need a few facts
about the moduli spaces. Fix a curve $C$ of genus $g$ and a group $G$.
They will be omitted from the notation.
\begin{lem} \label{properties}\mbox{}
\begin{enumerate}
\item $c_1(T_{\M}) = h\om$,
\item $p(T_{\M}) = c(\on{Ad} U_z)^{2(g-1)}=\prod_{\al\in\De} (1+\al)^{2(g-1)}$,
\item $\ahat(T_{\M}) =
\prod_{\al\in\De^+}\frwd{\al/2}{\sinh(\al/2)}^{2(g-1)}$.
\end{enumerate}
Here $h$ is the dual Coxeter number of $G$, $p$ denotes the total
Pontryagin class, $c$ the total Chern class,  $U_z$ is the bundle over
$\M$ obtained by restricting the universal principal bundle $U$ to a
slice $z\times\M$ for some $z\in C$  and  $\on{Ad} U_z$ is the vector bundle
associated to $U_z$ via the adjoint representation of $G$.
\end{lem}

For the proof of the first two statements in some partial cases see
\cite{AB}. The second statement
statement follows from the Kodaira-Spencer construction.  From the second
statement we find that the Pontryagin roots of $T_\M$ are the roots of the
Lie algebra $\g$, and this in turn implies the third statement.

Finally, we can put everything together. We will calculate $\dim
\cH^0(\M_C^G, \cL^k)$. Again, consider $G=\SL_3$ for simplicity.
First, the Kodaira vanishing theorem applies to
$\cL^k$, because the canonical bundle of $\M$ is negative, (this follows
from the \lemref{properties}(1), see also
\cite{BSz}). Thus we can replace the dimension of $\cH^0$ by the
Euler characteristic, and apply the Riemann-Roch theorem:
$$\dim\cH^0(\M_C^G, \cL^k)=
\chi(\M_C^G, \cL^k) = \int_{\M} e^{k\om} \Todd(\M). $$ According to
\lemref{properties}(1), and using the standard shifting trick we can rewrite
this integral as  $$\int_{\M} e^{(k+h)\om} \ahat(\M). $$

We can calculate this integral using \eqref{special} and
\lemref{properties}(3), and the result is exactly $V^G_k(g)$ according to
\eqref{verlslt}. This proves \eqref{main}.\qed

\end{document}